\documentclass[runningheads]{cl2emult}

\usepackage{makeidx}  
\usepackage{graphicx} 
\usepackage{subeqnar} 
\usepackage{multicol} 
\usepackage{cropmark} 
\usepackage{math}     
\makeindex            

\usepackage{amsmath,amsfonts}   
\def\GPunkt(#1,#2){\unskip
  \raise#2 \Einheit\hbox to0pt{\hskip#1 \Einheit
          \raise-2.5pt\hbox to0pt{\hss$\scriptstyle\bullet$\hss}\hss}}

\def\P{{\mathcal P}}
\def\O{{\mathcal O}}

\def\Z{{\mathbb Z}}
\def\PP{{\mathbb P}}
\def\R{{\mathbb R}}
\def\C{{\mathbb C}}
\def\rank{\operatorname{rank}}
\def\Char{\operatorname{char}}
\def\codim{\operatorname{codim}}
\def\GF{\operatorname{GF}}

\def\NE{\operatorname{NE}}
\def\Pf{\operatorname{Pf}}
\def\po#1#2{(#1)_#2}

\catcode`\@=11

\thinlines
\newskip\Einheit \Einheit=0.6cm
\newcount\xcoord \newcount\ycoord
\newdimen\xdim \newdimen\ydim \newdimen\PfadD@cke \newdimen\Pfadd@cke
\def\PfadDicke#1{\PfadD@cke#1 \divide\PfadD@cke by2 \Pfadd@cke\PfadD@cke \multiply\PfadD@cke by2}
\long\def\LOOP#1\REPEAT{\def\BODY{#1}\ITERATE}
\def\ITERATE{\BODY \let\next\ITERATE \else\let\next\relax\fi \next}
\let\REPEAT=\fi
\def\Punkt{\hbox{\raise-2pt\hbox to0pt{\hss\scriptsize$\bullet$\hss}}}
\def\DuennPunkt(#1,#2){\unskip
  \raise#2 \Einheit\hbox to0pt{\hskip#1 \Einheit
          \raise-2.5pt\hbox to0pt{\hss\normalsize$\bullet$\hss}\hss}}
\def\NormalPunkt(#1,#2){\unskip
  \raise#2 \Einheit\hbox to0pt{\hskip#1 \Einheit
          \raise-3pt\hbox to0pt{\hss\large$\bullet$\hss}\hss}}
\def\DickPunkt(#1,#2){\unskip
  \raise#2 \Einheit\hbox to0pt{\hskip#1 \Einheit
          \raise-4pt\hbox to0pt{\hss\Large$\bullet$\hss}\hss}}
\def\Kreis(#1,#2){\unskip
  \raise#2 \Einheit\hbox to0pt{\hskip#1 \Einheit
          \raise-4pt\hbox to0pt{\hss\Large$\circ$\hss}\hss}}
\def\Diagonale(#1,#2)#3{\unskip\leavevmode
  \xcoord#1\relax \ycoord#2\relax
      \raise\ycoord \Einheit\hbox to0pt{\hskip\xcoord \Einheit
         \unitlength\Einheit
         \line(1,1){#3}\hss}}
\def\AntiDiagonale(#1,#2)#3{\unskip\leavevmode
  \xcoord#1\relax \ycoord#2\relax \advance\xcoord by -0.05\relax
      \raise\ycoord \Einheit\hbox to0pt{\hskip\xcoord \Einheit
         \unitlength\Einheit
         \line(1,-1){#3}\hss}}
\def\Pfad(#1,#2),#3\endPfad{\unskip\leavevmode
  \xcoord#1 \ycoord#2 \thicklines\ZeichnePfad#3\endPfad\thinlines}
\def\ZeichnePfad#1{\ifx#1\endPfad\let\next\relax
  \else\let\next\ZeichnePfad
    \ifnum#1=1
      \raise\ycoord \Einheit\hbox to0pt{\hskip\xcoord \Einheit
         \vrule height\Pfadd@cke width1 \Einheit depth\Pfadd@cke\hss}%
      \advance\xcoord by 1
    \else\ifnum#1=2
      \raise\ycoord \Einheit\hbox to0pt{\hskip\xcoord \Einheit
        \hbox{\hskip-\PfadD@cke\vrule height1 \Einheit width\PfadD@cke depth0pt}\hss}%
      \advance\ycoord by 1
    \else\ifnum#1=3
      \raise\ycoord \Einheit\hbox to0pt{\hskip\xcoord \Einheit
         \unitlength\Einheit
         \line(1,1){1}\hss}
      \advance\xcoord by 1
      \advance\ycoord by 1
    \else\ifnum#1=4
      \raise\ycoord \Einheit\hbox to0pt{\hskip\xcoord \Einheit
         \unitlength\Einheit
         \line(1,-1){1}\hss}
      \advance\xcoord by 1
      \advance\ycoord by -1
    \fi\fi\fi\fi
  \fi\next}
\def\hSSchritt{\leavevmode\raise-.4pt\hbox to0pt{\hss.\hss}\hskip.2\Einheit
  \raise-.4pt\hbox to0pt{\hss.\hss}\hskip.2\Einheit
  \raise-.4pt\hbox to0pt{\hss.\hss}\hskip.2\Einheit
  \raise-.4pt\hbox to0pt{\hss.\hss}\hskip.2\Einheit
  \raise-.4pt\hbox to0pt{\hss.\hss}\hskip.2\Einheit}
\def\vSSchritt{\vbox{\baselineskip.2\Einheit\lineskiplimit0pt
\hbox{.}\hbox{.}\hbox{.}\hbox{.}\hbox{.}}}
\def\DSSchritt{\leavevmode\raise-.4pt\hbox to0pt{%
  \hbox to0pt{\hss.\hss}\hskip.2\Einheit
  \raise.2\Einheit\hbox to0pt{\hss.\hss}\hskip.2\Einheit
  \raise.4\Einheit\hbox to0pt{\hss.\hss}\hskip.2\Einheit
  \raise.6\Einheit\hbox to0pt{\hss.\hss}\hskip.2\Einheit
  \raise.8\Einheit\hbox to0pt{\hss.\hss}\hss}}
\def\dSSchritt{\leavevmode\raise-.4pt\hbox to0pt{%
  \hbox to0pt{\hss.\hss}\hskip.2\Einheit
  \raise-.2\Einheit\hbox to0pt{\hss.\hss}\hskip.2\Einheit
  \raise-.4\Einheit\hbox to0pt{\hss.\hss}\hskip.2\Einheit
  \raise-.6\Einheit\hbox to0pt{\hss.\hss}\hskip.2\Einheit
  \raise-.8\Einheit\hbox to0pt{\hss.\hss}\hss}}
\def\SPfad(#1,#2),#3\endSPfad{\unskip\leavevmode
  \xcoord#1 \ycoord#2 \ZeichneSPfad#3\endSPfad}
\def\ZeichneSPfad#1{\ifx#1\endSPfad\let\next\relax
  \else\let\next\ZeichneSPfad
    \ifnum#1=1
      \raise\ycoord \Einheit\hbox to0pt{\hskip\xcoord \Einheit
         \hSSchritt\hss}%
      \advance\xcoord by 1
    \else\ifnum#1=2
      \raise\ycoord \Einheit\hbox to0pt{\hskip\xcoord \Einheit
        \hbox{\hskip-2pt \vSSchritt}\hss}%
      \advance\ycoord by 1
    \else\ifnum#1=3
      \raise\ycoord \Einheit\hbox to0pt{\hskip\xcoord \Einheit
         \DSSchritt\hss}
      \advance\xcoord by 1
      \advance\ycoord by 1
    \else\ifnum#1=4
      \raise\ycoord \Einheit\hbox to0pt{\hskip\xcoord \Einheit
         \dSSchritt\hss}
      \advance\xcoord by 1
      \advance\ycoord by -1
    \fi\fi\fi\fi
  \fi\next}
\def\Koordinatenachsen(#1,#2){\unskip
 \hbox to0pt{\hskip-.5pt\vrule height#2 \Einheit width.5pt depth1 \Einheit}%
 \hbox to0pt{\hskip-1 \Einheit \xcoord#1 \advance\xcoord by1
    \vrule height0.25pt width\xcoord \Einheit depth0.25pt\hss}}
\def\Koordinatenachsen(#1,#2)(#3,#4){\unskip
 \hbox to0pt{\hskip-.5pt \ycoord-#4 \advance\ycoord by1
    \vrule height#2 \Einheit width.5pt depth\ycoord \Einheit}%
 \hbox to0pt{\hskip-1 \Einheit \hskip#3\Einheit 
    \xcoord#1 \advance\xcoord by1 \advance\xcoord by-#3 
    \vrule height0.25pt width\xcoord \Einheit depth0.25pt\hss}}
\def\Gitter(#1,#2){\unskip \xcoord0 \ycoord0 \leavevmode
  \LOOP\ifnum\ycoord<#2
    \loop\ifnum\xcoord<#1
      \raise\ycoord \Einheit\hbox to0pt{\hskip\xcoord \Einheit\Punkt\hss}%
      \advance\xcoord by1
    \repeat
    \xcoord0
    \advance\ycoord by1
  \REPEAT}
\def\Gitter(#1,#2)(#3,#4){\unskip \xcoord#3 \ycoord#4 \leavevmode
  \LOOP\ifnum\ycoord<#2
    \loop\ifnum\xcoord<#1
      \raise\ycoord \Einheit\hbox to0pt{\hskip\xcoord \Einheit\Punkt\hss}%
      \advance\xcoord by1
    \repeat
    \xcoord#3
    \advance\ycoord by1
  \REPEAT}
\def\Label#1#2(#3,#4){\unskip \xdim#3 \Einheit \ydim#4 \Einheit
  \def\lo{\advance\xdim by-.5 \Einheit \advance\ydim by.5 \Einheit}%
  \def\llo{\advance\xdim by-.25cm \advance\ydim by.5 \Einheit}%
  \def\loo{\advance\xdim by-.5 \Einheit \advance\ydim by.25cm}%
  \def\o{\advance\ydim by.25cm}%
  \def\ro{\advance\xdim by.5 \Einheit \advance\ydim by.5 \Einheit}%
  \def\rro{\advance\xdim by.25cm \advance\ydim by.5 \Einheit}%
  \def\roo{\advance\xdim by.5 \Einheit \advance\ydim by.25cm}%
  \def\l{\advance\xdim by-.30cm}%
  \def\r{\advance\xdim by.30cm}%
  \def\lu{\advance\xdim by-.5 \Einheit \advance\ydim by-.6 \Einheit}%
  \def\llu{\advance\xdim by-.25cm \advance\ydim by-.6 \Einheit}%
  \def\luu{\advance\xdim by-.5 \Einheit \advance\ydim by-.30cm}%
  \def\u{\advance\ydim by-.30cm}%
  \def\ru{\advance\xdim by.5 \Einheit \advance\ydim by-.6 \Einheit}%
  \def\rru{\advance\xdim by.25cm \advance\ydim by-.6 \Einheit}%
  \def\ruu{\advance\xdim by.5 \Einheit \advance\ydim by-.30cm}%
  #1\raise\ydim\hbox to0pt{\hskip\xdim
     \vbox to0pt{\vss\hbox to0pt{\hss$#2$\hss}\vss}\hss}%
}
\catcode`\@=12
\begin{document}
\title*{The Hilbert Series of Pfaffian Rings\thanks{The first author was partially supported by a 
`Career Award' grant from AICTE, New Delhi 
    and an IRCC grant from IIT Bombay. The second author was partially
supported by the Austrian
Science Foundation FWF, grant P13190-MAT.}}
\toctitle{The Hilbert series of Pfaffian rings}
%
%
\titlerunning{The Hilbert series of Pfaffian rings}
%
\author{Sudhir~R.~Ghorpade\inst{1}
\and Christian~Krattenthaler\inst{2}}
\authorrunning{S. R. Ghorpade and C. Krattenthaler}
%
%
\institute{
Department of Mathematics, \\
Indian Institute of Technology Bombay,\\
Powai, Mumbai 400076, India\\
e-mail: srg@math.iitb.ac.in\\
WWW: {\tt http://www.math.iitb.ac.in/\~{}srg}\\
\and 
Institut f\"ur Mathematik der Universit\"at Wien,\\
Strudlhofgasse 4, A-1090 Wien, Austria.\\
e-mail: KRATT@Ap.Univie.Ac.At\\
WWW: \tt http://www.mat.univie.ac.at/People/kratt
}

\maketitle              

\bigskip

\centerline{\it Dedicated to Professor Shreeram Abhyankar on his seventieth 
birthday}

\medskip

\begin{abstract}
We give three determinantal expressions for the Hilbert series
as well as the Hilbert function of a Pfaffian ring, and a closed form 
product formula for its multiplicity.
An appendix outlining some basic facts about degeneracy
loci and applications to multiplicity formulae for Pfaffian rings is
also included. 
\end{abstract}

\section{Introduction}
\label{Sec:Intro}

It is well known that the determinant of a skew-symmetric matrix of odd order is zero
whereas the determinant of a skew-symmetric matrix of even order is the square of a
polynomial in its entries, known as the {\em Pfaffian}. Combinatorially, a Pfaffian
may be described as the signed weight generating function of a complete graph.
We consider in this paper the
{\em Pfaffian ideals}, which are the ideals generated by the Pfaffians of a fixed size
in a generic skew-symmetric matrix, and the corresponding quotients of polynomial
rings, called {\em Pfaffian rings}. 
(See Section~\ref{Sec:Defns} for a more precise description.) 

Pfaffian rings have been studied by several authors and are known to possess a number
of nice properties. For example, Pfaffian rings are Cohen-Macaulay normal 
domains, which are, in fact,  factorial and Gorenstein (cf\@. \cite{AvraAA,KlLaAA,MariAA}).
Height, depth and in some cases, the minimal resolution of Pfaffian ideals is 
known (cf\@. \cite{JoPrAA,PragAA}). 
The singular 
locus of Pfaffian rings is known (cf\@. \cite{KlLaAA}) and the arithmetical
rank of Pfaffian ideals, i.e., the minimal number of equations needed 
to define the corresponding variety, is known as well (cf\@. \cite{BariAA}). 
Pfaffian rings arise in Invariant Theory as the ring of invariants of the
symplectic group (see, for example, \cite[Sec. 6]{DePrAA}).  In this 
connection, it is shown in \cite{DePrAA} (see also \cite[p. 53]{DeEiPr}) that 
there is a natural partial order on the set of Pfaffians of a skew-symmetric 
matrix and the corresponding polynomial ring is an ASL (algebra with 
straightening law) on the poset of Pfaffians.
The poset structure suggests the study of 
more general Pfaffian ideals (namely, those cogenerated by a Pfaffian), 
and algebraic properties of the corresponding residue class rings have
also been investigated (cf\@. \cite{BaetAA,DeNeAB}). 
Gr\"obner bases for Pfaffian rings 
have been constructed by Herzog and Trung \cite{HeTrAA}, who also derived
combinatorial formulae, in terms of the face numbers of an associated
simplicial complex, for the Hilbert function,
and a determinantal formula for the multiplicity.
An explicit expression for the Hilbert series of a Pfaffian ring has been
found by De Negri \cite[Theorem~3.5.1]{DeNeAA}, by establishing a
link between Pfaffian rings and ladder determinantal rings.

The purpose of our article is to record a few facts 
in this area that have been overlooked previously. 
First, we show that one can derive
an expression for the Hilbert series directly from the above
mentioned results of Herzog and Trung, if we combine them with
results from \cite{KratBE} and \cite{KrPrAA} 
on the enumeration of
nonintersecting lattice paths with respect to turns. 
The key in the derivation is to express the Hilbert
series in terms of a generating function for nonintersecting lattice
paths, in which the lattice paths are weighted by their number of
turns (see Proposition~\ref{Prop1}).
Second, in the process we not only recover De Negri's result, 
but also obtain 
two alternative expressions (see Theorem~\ref{MainThm1}). 
Third, we show that Herzog and Trung's determinantal formula
for the multiplicity actually simplifies to a nice closed product
formula (see Theorem~\ref{MainThm2}), thus answering a question raised in
\cite[p.~29]{HeTrAA}. This gives, for example, a formula for
the multiplicity of a generic Gorenstein ideal of codimension $3$ as a
trivial consequence.

Towards the completion of this work, we learned that formulae for the
multiplicity of Pfaffian rings can also be obtained by geometric methods
as a special case of the formulae for the fundamental classes of 
degeneracy loci of certain maps of vector bundles. 
The methods used in this paper are,
however, completely different, and also characteristic free. Moreover, the
resulting formulae are also somewhat distinct. Nevertheless, it seems
worthwhile to know the various formulae and the methods used to obtain them.
Accordingly, for the convenience of the reader, we have included a fairly 
self-contained appendix at the end of this paper in which
the basic ideas related to degeneracy loci and the resulting multiplicity
formulae are described. 

This paper is organized as follows. 
In the next section we review the definition of a Pfaffian ring, and
we introduce the lattice path notation that we are going to use
throughout the paper. The central part of the paper is 
Section~\ref{Sec:HilbSeries}, 
in which we establish the connection between the Hilbert series of a
Pfaffian ring and the enumeration of nonintersecting lattice
paths with a given number of turns. 
Section~\ref{Sec:MainResults} then contains our main results, the explicit formula
for the Hilbert series, and the closed form expression for the multiplicity. 
We close by discussing some applications and related work.

\section{Definitions}
\label{Sec:Defns}
Let $X=(X_{i,j})_{1\le i, j\le n}$ be a skew-symmetric $n\times n$ 
matrix where $\{X_{i,j}: i< j\}$ are independent
indeterminates over a field $K$. 
Let $K[X]$ denote the ring of all polynomials 
in the $X_{i,j}$'s, 
with coefficients in $K$, and let $I_{r+1}(X)$
be the ideal of $K[X]$ that is generated by all $(2r+2)\times (2r+2)$
Pfaffian minors\footnote{The {\it Pfaffian} $\Pf(A)$ of a
skew-symmetric $(2m)\times(2m)$ matrix $A$ is defined by
$$\Pf(A)=\sum _{\pi} ^{}(-1)^{c(\pi)}\prod _{(ij)\in \pi} ^{}A_{ij},$$
where the sum is over all perfect matchings $\pi$ of the complete
graph on $2m$ vertices, where $c(\pi)$ is the {\it crossing
number} of $\pi$, and where the product is over all edges $(ij)$,
$i<j$, in the matching $\pi$ (see e.g., 
\cite[Sec.~2]{StemAE}). A {\it Pfaffian minor} of a skew-symmetric matrix $X$ 
is the Pfaffian of a submatrix of $X$ consisting of the rows and
columns indexed by $i_1$, $i_2$, \dots, $i_{2r+2}$, for some
$i_1<i_2<\dots<i_{2r+2}$.} of $X$.
The ideal $I_{r+1}(X)$ is called a {\it Pfaffian
ideal}. The associated {\it Pfaffian
ring\/} is $R_{r+1}(X):=K[X]/I_{r+1}(X)$. 

Throughout the paper by a {\it lattice path} we mean a lattice path 
in the plane integer lattice $\Z^2$ ($\Z$ denoting the set of
integers) consisting of unit horizontal and vertical steps
in the positive direction. In the sequel we shall frequently
refer to them as {\it paths}. See Figure~1 for an example of a path
$P_0$ from $(1,-1)$ to $(6,6)$.
We shall frequently abbreviate the fact that a path $P$ goes from
$A$ to $E$ by writing $P:A\to E$.

\begin{figure}[h]
$$\Koordinatenachsen(8,8)(0,0)
\Gitter(8,8)(0,0)
\Pfad(1,-1),221221112122\endPfad
\DickPunkt(1,-1)
\DickPunkt(6,6)
\Label\ro{P_0}(3,3)
\hskip4cm
$$
\caption{}
\end{figure}

Also, given lattice points $A$ and $E$, we denote the set of all lattice
paths from $A$ to $E$ by $\P(A\to E)$.   
A family $(P_1,P_2,\dots,P_r)$ of lattice paths is said to be 
{\it nonintersecting\/} if no
two lattice paths of this family have a point in common. 

A point in a path $P$ which is the end point of a vertical
step and at the same time
 the starting point of a horizontal step will be called
a {\it North-East turn} ({\it NE-turn} for short) of the path $P$. 
The NE-turns of
the path in Figure~1 are $(1,1)$, $(2,3)$, and $(5,4)$. 
We write $\NE(P)$ for the number of NE-turns of $P$. Also, given a family
$\mathbf P=(P_1,P_2,\dots,P_n)$ of paths $P_i$, we write $\NE(\mathbf P)$ for the
number $\sum _{i=1} ^{n}\NE(P_i)$ of all NE-turns in the family.
Finally, given any weight function $w$ defined on a set $\mathcal M$ 
(and taking values in a ring),
by the generating function $\GF(\mathcal M;w)$ we mean $\sum
_{x\in\mathcal M} ^{}w(x)$.

\section{The Hilbert series of a Pfaffian ring and nonintersecting lattice paths}
\label{Sec:HilbSeries}
In this section we establish the central result of this paper, the
connection between the Hilbert series of a
Pfaffian ring and enumeration of nonintersecting lattice
paths with a given number of NE turns.
\begin{proposition}
\label{Prop1}
Let $A_i=(r+i-1,r-i+1)$ and $E_i=(n-r+i-1,n-r-i+1)$, $i=1,2,\dots,r$, be
lattice points.
Then the Hilbert series of the Pfaffian
ring $R_{r+1}(X)=K[X]/I_{r+1}(X)$ equals
\begin{equation} \label{e3.1} 
\sum _{\ell=0} ^{\infty}\dim_K R_{r+1}(X)_\ell\,z^\ell
=\frac {\GF(\P^+(\mathbf A\to \mathbf E);z^{\NE(.)})}
{(1-z)^{r(2n-2r-1)}},
\end{equation}
where $R_{r+1}(X)_\ell$ denotes the homogeneous component of degree
$\ell$ in $R_{r+1}(X)$, and where $\P^+(\mathbf A\to\mathbf E)$ denotes
the set of all 
families $(P_1,P_2,\dots,P_r)$ of nonintersecting lattice
paths, where $P_i$ runs from $A_i$ to $E_i$ and
never passes above the diagonal $x=y$, $i=1,2,\dots,r$.
\end{proposition}
\begin{proof}
We use some results of Herzog and Trung
\cite{HeTrAA}. Our arguments are completely parallel to the
arguments in the second proof of Theorem~2 in \cite{KrPrAA}.

In Section~5 of \cite{HeTrAA}, Pfaffian
rings are introduced and investigated. 
It is shown there
that for a suitable term order (order on monomials), the ideal $I_{r+1}(X)^*$
of leading monomials of $I_{r+1}(X)$ is generated by
square-free monomials. Thus $K[X]/I_{r+1}(X)^*$ may be viewed as a
Stanley-Reisner ring of a certain simplicial complex $\Delta_{r+1}$. 
The faces
of this simplicial complex $\Delta_{r+1}$ are described in Lemma~5.3 of
\cite{HeTrAA}. Namely, translated into a less formal language, the
faces are sets $S$ of integer lattice points in the (upper) triangular region 
$\{(x,y):1\le x<y\le n\}$, such that 
\begin{equation} \label{e3.2} 
\begin{matrix} 
\vbox{\hsize9cm\noindent
a sequence 
$(i_1,j_1)$, $(i_2,j_2)$, \dots, $(i_k,j_k)$ of
elements of $S$ with $i_1<i_2<\dots<i_k$ and $j_1>j_2>\dots>j_k$ does
not contain more than $r$ elements.
}\end{matrix}
\end{equation}

\begin{figure}[h]
$$
\Einheit.4cm
\Koordinatenachsen(12,12)(0,0)
\Kreis(0,0)
\GPunkt(1,2)
\GPunkt(1,3)
\GPunkt(1,4)
\GPunkt(1,5)
\GPunkt(1,6)
\GPunkt(1,7)
\GPunkt(1,8)
\GPunkt(1,9)
\GPunkt(1,10)
\GPunkt(1,11)
\GPunkt(1,12)
\GPunkt(2,3)
\GPunkt(2,4)
\GPunkt(2,5)
\GPunkt(2,6)
\GPunkt(2,7)
\GPunkt(2,8)
\GPunkt(2,9)
\GPunkt(2,10)
\GPunkt(2,11)
\GPunkt(2,12)
\GPunkt(3,4)
\GPunkt(3,5)
\GPunkt(3,6)
\GPunkt(3,7)
\GPunkt(3,8)
\GPunkt(3,9)
\GPunkt(3,10)
\GPunkt(3,11)
\GPunkt(3,12)
\GPunkt(4,5)
\GPunkt(4,6)
\GPunkt(4,7)
\GPunkt(4,8)
\GPunkt(4,9)
\GPunkt(4,10)
\GPunkt(4,11)
\GPunkt(4,12)
\GPunkt(5,6)
\GPunkt(5,7)
\GPunkt(5,8)
\GPunkt(5,9)
\GPunkt(5,10)
\GPunkt(5,11)
\GPunkt(5,12)
\GPunkt(6,7)
\GPunkt(6,8)
\GPunkt(6,9)
\GPunkt(6,10)
\GPunkt(6,11)
\GPunkt(6,12)
\GPunkt(7,8)
\GPunkt(7,9)
\GPunkt(7,10)
\GPunkt(7,11)
\GPunkt(7,12)
\GPunkt(8,9)
\GPunkt(8,10)
\GPunkt(8,11)
\GPunkt(8,12)
\GPunkt(9,10)
\GPunkt(9,11)
\GPunkt(9,12)
\GPunkt(10,11)
\GPunkt(10,12)
\GPunkt(11,12)
\DickPunkt(1,3)
\DickPunkt(2,4)
\DickPunkt(3,4)
\DickPunkt(4,5)
\DickPunkt(1,6)
\DickPunkt(1,7)
\DickPunkt(3,7)
\DickPunkt(5,7)
\DickPunkt(5,8)
\DickPunkt(2,9)
\DickPunkt(3,9)
\DickPunkt(7,9)
\DickPunkt(3,10)
\DickPunkt(6,10)
\DickPunkt(2,11)
\DickPunkt(3,11)
\DickPunkt(8,11)
\DickPunkt(10,11)
\DickPunkt(5,12)
\DickPunkt(9,12)
\hskip5cm
$$
\caption{}
\end{figure}

An example of such a point set $S$, with $n=12$ and $r=3$, is the set
\begin{multline*}
\{
(1,3),(2,4),(3,4),(4,5),(1,6),(1,7),(3,7),(5,7),(5,8),(2,9),
(3,9),\\(7,9),
(3,10),(6,10),(2,11),(3,11),(8,11),(10,11),(5,12),(9,12)
\}.
\end{multline*}
A geometric realization of this set is contained in Figure~2,
the elements of the set $S$ being indicated by bold dots, 
the small dots indicating the triangular region of which $S$ is a
subset.

As usual, let $f_i$ be the number of $i$-dimensional faces of $\Delta_{r+1}$,
i.e., the number of such sets $S$ of cardinality $i+1$. 
Corollary~5.2 of \cite{HeTrAA} says that the dimensions of the
homogeneous components of the Pfaffian ring $R_{r+1}(X)$ can be
expressed in terms of the face numbers $f_i$, namely there holds 
$$
\dim_K R_{r+1}(X)_\ell=\sum _{i\ge0} ^{}\binom{\ell-1}i f_i 
\quad {\rm for \ every \ } \ell \ge 0.
$$

Now consider such a set $S$. For convenience we apply to it the
mapping $(x,y)\to (y-1,x)$ (i.e., the reflection in the main diagonal
followed by a shift by $1$ in the negative $x$-direction). 
Thus we obtain a point set, $\tilde S$ say, in the (lower) triangular region
$\{1\le y\le x\le n-1\}$. See
Figure~3.a for the result when this mapping is applied to our example
point set in Figure~2.

\begin{figure}[h]
$$
\Einheit.4cm
\Koordinatenachsen(12,12)(0,0)
\Kreis(0,0)
\GPunkt(1,1)
\GPunkt(2,1)
\GPunkt(3,1)
\GPunkt(4,1)
\GPunkt(5,1)
\GPunkt(6,1)
\GPunkt(7,1)
\GPunkt(8,1)
\GPunkt(9,1)
\GPunkt(10,1)
\GPunkt(11,1)
\GPunkt(2,2)
\GPunkt(3,2)
\GPunkt(4,2)
\GPunkt(5,2)
\GPunkt(6,2)
\GPunkt(7,2)
\GPunkt(8,2)
\GPunkt(9,2)
\GPunkt(10,2)
\GPunkt(11,2)
\GPunkt(3,3)
\GPunkt(4,3)
\GPunkt(5,3)
\GPunkt(6,3)
\GPunkt(7,3)
\GPunkt(8,3)
\GPunkt(9,3)
\GPunkt(10,3)
\GPunkt(11,3)
\GPunkt(4,4)
\GPunkt(5,4)
\GPunkt(6,4)
\GPunkt(7,4)
\GPunkt(8,4)
\GPunkt(9,4)
\GPunkt(10,4)
\GPunkt(11,4)
\GPunkt(5,5)
\GPunkt(6,5)
\GPunkt(7,5)
\GPunkt(8,5)
\GPunkt(9,5)
\GPunkt(10,5)
\GPunkt(11,5)
\GPunkt(6,6)
\GPunkt(7,6)
\GPunkt(8,6)
\GPunkt(9,6)
\GPunkt(10,6)
\GPunkt(11,6)
\GPunkt(7,7)
\GPunkt(8,7)
\GPunkt(9,7)
\GPunkt(10,7)
\GPunkt(11,7)
\GPunkt(8,8)
\GPunkt(9,8)
\GPunkt(10,8)
\GPunkt(11,8)
\GPunkt(9,9)
\GPunkt(10,9)
\GPunkt(11,9)
\GPunkt(10,10)
\GPunkt(11,10)
\GPunkt(11,11)
\DickPunkt(2,1)
\DickPunkt(3,2)
\DickPunkt(3,3)
\DickPunkt(4,4)
\DickPunkt(5,1)
\DickPunkt(6,1)
\DickPunkt(6,3)
\DickPunkt(6,5)
\DickPunkt(7,5)
\DickPunkt(8,2)
\DickPunkt(8,3)
\DickPunkt(8,7)
\DickPunkt(9,3)
\DickPunkt(9,6)
\DickPunkt(10,2)
\DickPunkt(10,3)
\DickPunkt(10,8)
\DickPunkt(10,10)
\DickPunkt(11,5)
\DickPunkt(11,9)
\hbox{\hskip5cm}
$$
\centerline{\small a. The set of points reflected}
\vskip6pt
$$
\Koordinatenachsen(12,12)(0,0)
\Kreis(0,0)
\GPunkt(1,1)
\GPunkt(2,1)
\GPunkt(3,1)
\GPunkt(4,1)
\GPunkt(5,1)
\GPunkt(6,1)
\GPunkt(7,1)
\GPunkt(8,1)
\GPunkt(9,1)
\GPunkt(10,1)
\GPunkt(11,1)
\GPunkt(2,2)
\GPunkt(3,2)
\GPunkt(4,2)
\GPunkt(5,2)
\GPunkt(6,2)
\GPunkt(7,2)
\GPunkt(8,2)
\GPunkt(9,2)
\GPunkt(10,2)
\GPunkt(11,2)
\GPunkt(3,3)
\GPunkt(4,3)
\GPunkt(5,3)
\GPunkt(6,3)
\GPunkt(7,3)
\GPunkt(8,3)
\GPunkt(9,3)
\GPunkt(10,3)
\GPunkt(11,3)
\GPunkt(4,4)
\GPunkt(5,4)
\GPunkt(6,4)
\GPunkt(7,4)
\GPunkt(8,4)
\GPunkt(9,4)
\GPunkt(10,4)
\GPunkt(11,4)
\GPunkt(5,5)
\GPunkt(6,5)
\GPunkt(7,5)
\GPunkt(8,5)
\GPunkt(9,5)
\GPunkt(10,5)
\GPunkt(11,5)
\GPunkt(6,6)
\GPunkt(7,6)
\GPunkt(8,6)
\GPunkt(9,6)
\GPunkt(10,6)
\GPunkt(11,6)
\GPunkt(7,7)
\GPunkt(8,7)
\GPunkt(9,7)
\GPunkt(10,7)
\GPunkt(11,7)
\GPunkt(8,8)
\GPunkt(9,8)
\GPunkt(10,8)
\GPunkt(11,8)
\GPunkt(9,9)
\GPunkt(10,9)
\GPunkt(11,9)
\GPunkt(10,10)
\GPunkt(11,10)
\GPunkt(11,11)
\DickPunkt(2,1)
\DickPunkt(3,2)
\DickPunkt(3,3)
\DickPunkt(4,4)
\DickPunkt(5,1)
\DickPunkt(6,1)
\DickPunkt(6,3)
\DickPunkt(6,5)
\DickPunkt(7,5)
\DickPunkt(8,2)
\DickPunkt(8,3)
\DickPunkt(8,7)
\DickPunkt(9,3)
\DickPunkt(9,6)
\DickPunkt(10,2)
\DickPunkt(10,3)
\DickPunkt(10,8)
\DickPunkt(10,10)
\DickPunkt(11,5)
\DickPunkt(11,9)
\Diagonale(-1,-1){13}
\Pfad(3,3),121121122122\endPfad
\Pfad(4,2),112111222122\endPfad
\Pfad(5,1),111211212222\endPfad
\SPfad(-1,7),4444444\endSPfad
\SPfad(6,12),444444\endSPfad
\Kreis(3,3)
\Kreis(4,2)
\Kreis(5,1)
\Kreis(9,9)
\Kreis(10,8)
\Kreis(11,7)
\Label\l{P_1\kern-5pt}(8,6)
\Label\l{P_2\kern-5pt}(9,4)
\Label\ro{P_3}(11,3)
\Label\l{A_1}(3,3)
\Label\lu{A_2}(4,2)
\Label\u{A_3}(5,1)
\Label\o{E_1}(9,9)
\Label\ro{E_2}(10,8)
\Label\r{E_3}(11,7)
\Label\ro{x=y\kern5pt}(12,12)
\Label\lo{x+y=2r}(-1,7)
\Label\lo{x+y=2n-2r}(6,12)
\hskip6.5cm
$$
\centerline{\small b. 
The corresponding family of nonintersecting lattice paths}
\vskip6pt
\caption{}
\end{figure}

Next we apply to $\tilde S$ a variant of Viennot's ``light and shadow procedure"
(see \cite{GhorAB,HeTrAA,SagaAL,VienAA}). 
This variant defines, for each such point set $\tilde S$ (and, thus,
for each point set $S$ in the upper triangular region), a
family $(P_1,P_2,\dots,P_r)$ of $r$ nonintersecting lattice paths,
$P_i$ running from $(r+i-1,r-i+1)$ to $(n-r+i-1,n-r-i+1)$,
$i=1,2,\dots,r$, in the following way. 

First we ignore everything
which is in the left-bottom corner of the triangular region to the
left/below of the line $x+y=2r$, and everything 
which is in the right-top corner of the triangular region to the
right/above of the line $x+y=2n-2r$. In our running example, these
two lines are indicated as dotted lines in Figure~3.b. 

Next we suppose that there is a light source being located in the
top-left corner. The {\it shadow} of a point
$(x,y)$ is defined to be the set of points $(x',y')\in\R^2$ ($\R$
denoting the set of real numbers) with $x\le x'$ and $y'\le y$. 
We consider the (top-left) {\it border} of the union of
the shadows of all the points of the set $S$ that are located inside the
strip between the two lines $x+y=2r$ and $x+y=2n-2r$. We also include
the shadows of the points $A_1=(r,r)$ and $E_1=(n-r,n-r)$. This border is a
lattice path, $P_1$ say, from $A_1$ to $E_1$. Now we remove all the points of the
set that lie on this path. Then the
light and shadow procedure is repeated with the remaining points. 
In the second run we also include the shadows of $A_2=(r+1,r-1)$
and $E_2=(n-r+1,n-r-1)$, etc. We stop after a total of $r$ iterations.
Thus we obtain exactly $r$ lattice paths, the $i$-th path, $P_i$ say, 
running from $A_i$ to $E_i$. Clearly, by construction, the paths have
the property that they are nonintersecting and that they never pass
above the diagonal $x=y$. In addition,
a moment's thought shows that condition \eqref{e3.2} guarantees that after
these $r$ iterations all the points of the set $S$ are exhausted.
Figure~3.b displays the lattice paths
which in our example are obtained by this procedure.

On the other hand, if we are given a family 
$(P_1,P_2,\dots,P_r)$ of $r$ nonintersecting
lattice paths, $P_i$ running from $A_i$ to $E_i$, 
$i=1,2,\dots,r$, with a total number of
exactly $m$ NE-turns, how many point sets $S$ of cardinality $i+1$ 
in the (upper) triangular region 
$\{(x,y):1\le x<y\le n\}$ satisfying \eqref{e3.2} are there which, after the
transformation $(x,y)\to (y-1,x)$ and subsequent light and shadow as
described above, generate the given family of nonintersecting lattice
paths? Clearly, every NE-turn of a path of the family must be
occupied by a point of $\tilde S$. Aside from that, any point on any of the $r$ paths, 
any point in the bottom-left corner cut off by $x+y=2r$, and any
point in the top-right corner cut off by $x+y=2n-2r$ may or may not be in
$\tilde S$.
Hence, if we denote by $d$ the total number of points in the union of
the $r$ paths and these two corner regions, then
there are exactly $\binom{d-m}{i+1-m}$
sets $S$ of cardinality $i+1$ that reduce to
$(P_1,P_2,\dots,P_r)$ under light and shadow. As an easy computation
shows, we have $d=2r^2+r(2n-4r-1)=r(2n-2r-1)$.

Hence, if $h_m$ denotes the number of all 
families $(P_1,P_2,\dots,P_r)$ of $r$ nonintersecting
lattice paths, $P_i$ running from $A_i$ to $E_i$, 
$i=1,2,\dots,r$, with a total number of
exactly $m$ NE-turns, we see that the Hilbert series equals
\begin{align*} \sum _{\ell=0} ^{\infty}\dim_K R_{r+1}(X)_\ell\, z^\ell&=
\sum _{\ell=0} ^{\infty}\bigg(\sum _{i\ge0}\binom{\ell-1}i f_i\bigg)z^\ell\\
&=\sum _{\ell=0} ^{\infty}\sum _{i\ge0}\binom{\ell-1}i
\bigg(\sum _{m=0} ^{i+1}\binom{d-m}{i+1-m}h_m\bigg)z^\ell\\
&=\sum _{m=0} ^{\infty}h_m\sum _{\ell=0}
^{\infty}z^\ell\sum _{i\ge0} ^{}\binom{\ell-1}i\binom{d-m}{d-i-1},
\end{align*}
and if we sum the inner sum by means of the Chu--Vandermonde summation 
(see e.g\@. \cite[Sec.~5.1, (5.27)]{GrKPAA}), then we obtain
$$
\sum _{\ell=0} ^{\infty}\dim_K R_{r+1}(X)_\ell\, z^\ell=
\sum _{m=0} ^{\infty}h_m\sum _{\ell=0}
^{\infty}z^\ell\binom{d+\ell-m-1}{d-1}.$$
In the sum over $\ell$, the terms for $\ell <m$ vanish, so that we
may sum over $\ell\ge m$. Application of the binomial theorem then
yields
$$\sum _{\ell=0} ^{\infty}\dim_K R_{r+1}(X)_\ell\, z^\ell=
\frac {\sum _{m=0} ^{\infty}h_mz^m} {(1-z)^d}.
$$
This is exactly \eqref{e3.1}.
\quad \quad \qed

\end{proof}

\section{The main results}
\label{Sec:MainResults}
Our determinantal formulae for the Hilbert series of a Pfaffian ring
are the following.
\begin{theorem}
\label{MainThm1}
The Hilbert series of the Pfaffian
ring $R_{r+1}(X)=\break K[X]/I_{r+1}(X)$ equals
\begin{multline} \label{e4.1} 
\sum _{\ell=0} ^{\infty}\dim_K R_{r+1}(X)_\ell\,z^\ell\\
=\frac {\det\limits_{1\le i,j\le r}\Big(\sum\limits _{k} ^{}\Big(
\binom {n-2r}{k+i-j}\binom {n-2r}{k}-
\binom {n-2r-1}{k-j}\binom {n-2r+1}{k+i}\Big)z^k\Big)}
{(1-z)^{r(2n-2r-1)}},
\end{multline}
or, alternatively,
\begin{multline} \label{e4.1a} 
\sum _{\ell=0} ^{\infty}\dim_K R_{r+1}(X)_\ell\,z^\ell\\
=\frac {\det\limits_{1\le i,j\le r}\Big(\sum\limits _{k} ^{}\Big(
\binom {n-2r+i-1}{k+i-j}\binom {n-2r+j-1}{k}-
\binom {n-2r-1}{k-j}\binom {n-2r+i+j-1}{k+i}\Big)z^k\Big)}
{(1-z)^{r(2n-2r-1)}},
\end{multline}
or, alternatively,
\begin{multline} \label{e4.2} 
\sum _{\ell=0} ^{\infty}\dim_K R_{r+1}(X)_\ell\,z^\ell\\
=\frac {z^{-\binom r2}
\det\limits_{1\le i,j\le r}\Big(\sum\limits _{k} ^{}\Big(
\binom {n-2r+i-1}{k}\binom {n-2r+j-1}{k}-
\binom {n-2r+i+j-3}{k-1}\binom {n-2r+1}{k+1}\Big)z^k\Big)}
{(1-z)^{r(2n-2r-1)}},\kern-5pt
\end{multline}
where, once more, $R_{r+1}(X)_\ell$ denotes the homogeneous component of degree
$\ell$ of $R_{r+1}(X)$.
\end{theorem}

\begin{proof}
In view of Proposition~\ref{Prop1}, we only have to solve the
problem of enumeration, with respect to NE-turns, 
of nonintersecting lattice paths that are bounded by a
diagonal line. This has been previously accomplished
in \cite{KratBE} and in \cite{KrPrAA}. To be precise, to show that
the generating
function $\GF(\P^+(\mathbf A\to \mathbf E);z^{\NE(.)})$ in the numerator
on the right-hand side of \eqref{e3.1} can be expressed by the determinant
on the right-hand side of \eqref{e4.1}, one sets $A_1^{(i)}=r+i-1$,
$A_2^{(i)}=r-i+1$, $E_1^{(i)}=n-r+i-1$, $E_2^{(i)}=n-r-i+1$ in
Theorem~2 of \cite{KratBE}, then multiplies the resulting expression
by $z^{K}$, and sums over all $K$. 

\begin{figure}[h]
$$
\Einheit.4cm
\Koordinatenachsen(12,12)(0,0)
\Kreis(0,0)
\GPunkt(1,1)
\GPunkt(2,1)
\GPunkt(3,1)
\GPunkt(4,1)
\GPunkt(5,1)
\GPunkt(6,1)
\GPunkt(7,1)
\GPunkt(8,1)
\GPunkt(9,1)
\GPunkt(10,1)
\GPunkt(11,1)
\GPunkt(2,2)
\GPunkt(3,2)
\GPunkt(4,2)
\GPunkt(5,2)
\GPunkt(6,2)
\GPunkt(7,2)
\GPunkt(8,2)
\GPunkt(9,2)
\GPunkt(10,2)
\GPunkt(11,2)
\GPunkt(3,3)
\GPunkt(4,3)
\GPunkt(5,3)
\GPunkt(6,3)
\GPunkt(7,3)
\GPunkt(8,3)
\GPunkt(9,3)
\GPunkt(10,3)
\GPunkt(11,3)
\GPunkt(4,4)
\GPunkt(5,4)
\GPunkt(6,4)
\GPunkt(7,4)
\GPunkt(8,4)
\GPunkt(9,4)
\GPunkt(10,4)
\GPunkt(11,4)
\GPunkt(5,5)
\GPunkt(6,5)
\GPunkt(7,5)
\GPunkt(8,5)
\GPunkt(9,5)
\GPunkt(10,5)
\GPunkt(11,5)
\GPunkt(6,6)
\GPunkt(7,6)
\GPunkt(8,6)
\GPunkt(9,6)
\GPunkt(10,6)
\GPunkt(11,6)
\GPunkt(7,7)
\GPunkt(8,7)
\GPunkt(9,7)
\GPunkt(10,7)
\GPunkt(11,7)
\GPunkt(8,8)
\GPunkt(9,8)
\GPunkt(10,8)
\GPunkt(11,8)
\GPunkt(9,9)
\GPunkt(10,9)
\GPunkt(11,9)
\GPunkt(10,10)
\GPunkt(11,10)
\GPunkt(11,11)
\Diagonale(-1,-1){13}
\Pfad(3,3),121121122122\endPfad
\Pfad(3,2),11121112221222\endPfad
\Pfad(3,1),1111121121222222\endPfad
\Kreis(3,3)
\Kreis(3,2)
\Kreis(3,1)
\Kreis(9,9)
\Kreis(10,9)
\Kreis(11,9)
\Label\l{A'_1}(3,3)
\Label\l{A'_2}(3,2)
\Label\l{A'_3}(3,1)
\Label\o{E'_1}(9,9)
\Label\o{E'_2}(10,9)
\Label\o{E'_3}(11,9)
\Label\ro{x=y\kern5pt}(12,12)
\hskip4.8cm
$$
\caption{}
\end{figure}

To show that it can be expressed
by the numerator on the right-hand side of \eqref{e4.1a}, 
respectively
of \eqref{e4.2}, 
we prepend $(i-1)$ horizontal
steps and append $(i-1)$ vertical steps to $P_i$.
Then, out of a
family of nonintersecting paths as in the statement of Proposition~\ref{Prop1}, 
we obtain a family
$(P'_1,P'_2,\dots,P'_r)$ of nonintersecting lattice
paths, where $P'_i$ runs from $A'_i=(r,r-i+1)$ to
$E'_i=(n-r+i-1,n-r)$ and does not pass above $x=y$, 
$i=1,2,\dots,r$. See Figure~4 for the
corresponding path family which is obtained out of the one in
Figure~3.b. Clearly, the number of the
latter families is exactly the same as the number of the former,
because the prepended and appended portions are ``forced," i.e., if 
$(P'_1,P'_2,\dots,P'_r)$ are nonintersecting, then they must
contain these prepended and appended portions.
Now one can either again apply Theorem~2 in \cite{KratBE}, this time with 
$A_1^{(i)}=r$,
$A_2^{(i)}=r-i+1$, $E_1^{(i)}=n-r+i-1$, $E_2^{(i)}=n-r$, multiply the
resulting expression by $z^{K}$, sum over all $K$, and thus obtain the
numerator in \eqref{e4.1a}, or apply
Theorem~2 in \cite{KrPrAA} in conjunction with Proposition~6, (6.6), in
\cite{KrPrAA} with $D=0$, and thus obtain the numerator in
\eqref{e4.2}.
\quad \quad \qed
\end{proof}

\medskip
\noindent
{\it Remark}. Formula \eqref{e4.2} had been found earlier by De Negri
\cite[Theorem~3.5.1]{DeNeAA}.
\medskip

For convenience, we make the resulting expressions for the Hilbert
{\it function} explicit.

\begin{corollary}
\label{Cor1}
The Hilbert function of the Pfaffian
ring $R_{r+1}(X)=\break K[X]/I_{r+1}(X)$ is given by 
\begin{equation} 
\label{e4.3}
\dim_K R_{r+1}(X)_\ell
=\sum_{k} F_k \binom {\ell+r(2n-2r-1)-k-1}{r(2n-2r-1)-1}
\end{equation} 
where for $k\in \Z$, the coefficient $F_k$ equals
$$
\sum_{k_1+\cdots+k_r = k }
\det\limits_{1\le i,j\le r}\bigg(
\binom {n-2r}{k_i+i-j}\binom {n-2r}{k_i}-
\binom {n-2r-1}{k_i-j}\binom {n-2r+1}{k_i+i}\bigg),
$$
or, alternatively,
\begin{equation} 
\label{e4.3a}
\dim_K R_{r+1}(X)_\ell
=\sum_{k} G_k \binom {\ell+r(2n-2r-1)-k-1}{r(2n-2r-1)-1}
\end{equation} 
where for $k\in \Z$, the coefficient $G_k$ equals
\begin{multline*}
\sum_{k_1+\cdots+k_r = k }
\det\limits_{1\le i,j\le r}\bigg(
\binom {n-2r+i-1}{k_i+i-j}\binom {n-2r+j-1}{k_i}
\\-
\binom {n-2r-1}{k_i-j}\binom {n-2r+i+j-1}{k_i+i}\bigg),
\end{multline*}
or, alternatively,
\begin{equation} 
\label{e4.4}
\dim_K R_{r+1}(X)_\ell
= \sum _{k} H_k \binom {\ell+r(2n-2r-1)+\binom r2- k -1}{r(2n-2r-1)-1}
\end{equation} 
where for $k\in \Z$, the coefficient $H_k$ equals
\begin{multline*} 
\sum_{k_1+\cdots+k_r = k }
\det\limits_{1\le i,j\le r}\bigg(
\binom {n-2r+i-1}{k_i}\binom {n-2r+j-1}{k_i} \\
\kern2cm-
 \binom {n-2r+i+j-3}{k_i-1}\binom {n-2r+1}{k_i+1}\bigg).
\end{multline*}
\end{corollary}

\medskip
\noindent
{\it Remarks}. (1) The sums over 
$k_1,k_2,\dots,k_r$ appearing in the corollary above
are in fact finite sums, because each of the $k_i$'s is bounded
above and below due to the binomial coefficients which appear in
the determinants. This shows also that $F_k$ as well as $G_k$ and
$H_k$ are
zero for all except finitely many $k\in \Z$, and consequently, the sums over $k$ in
\eqref{e4.3}, \eqref{e4.3a} and \eqref{e4.4} are also finite. 
Thus, in particular, the expressions \eqref{e4.3}, \eqref{e4.3a}
and \eqref{e4.4}
exhibit transparently that the Hilbert function is a polynomial in
$\ell$ for {\it all\/} $\ell$. 
This proves that the ideal $I_{r+1}(X)$ is {\it Hilbertian} in the sense of
Abhyankar \cite{AbhyAB}.
\medskip

(2) It may be interesting to note that a formula for the Hilbert function of $R_{r+1}(X)$ is already
known in the special case of $r=1$. Indeed, the ideal $I_2(X)$ of $4\times 4$ Pfaffians
in a $n\times n$ skew-symmetric matrix precisely equals the ideal of the Pl\"ucker
relations in the Grassmannian $G_{2,n}$ of $2$-planes in $n$-space (over $K$). The
Hilbert function of an arbitrary Grassmannian $G_{d,n}$ and, more generally, of any 
Schubert variety $\Omega_{\alpha}$ in $G_{d,n}$ was determined by Hodge \cite{HodgAA}
in 1943 (see also \cite{GhorAC}).  Using Hodge's formula in this special case
(e.g., putting $d=2$ and $\alpha_i = n-d+i$ in \cite[Theorem 6]{GhorAC}), we see
that $\dim_K R_{2}(X)_\ell$ equals
\begin{equation} 
\label{Eq:Ho1}
{\binom{ \ell + n - 2}{\ell} }^2 - 
{\binom{ \ell + n - 2}{\ell - 1} } {\binom{ \ell + n - 2}{\ell + 1} }.
\end{equation} 
On the other hand, in the case of $r=1$, the formula \eqref{e4.3} of Corollary~\ref{Cor1} 
reduces to the following seemingly more complicated expression:
\begin{equation} 
\label{Eq:Ho2}
\sum_{k} 
{\binom{2n + \ell - k - 4}{\ell -k} } \left[
{\binom{n-2}{k}}^2  -  {\binom{n-3}{k-1}} {\binom{n-1}{k+1}} \right].
\end{equation} 
The resulting identity of \eqref{Eq:Ho1} and \eqref{Eq:Ho2} is not difficult to verify
directly. In fact, both \eqref{Eq:Ho1} and \eqref{Eq:Ho2} are differences of two terms 
and the corresponding terms are also equal to each other. Indeed,
using the standard hypergeometric notation
\begin{equation*} 
{}_r F_s\!\left[\begin{matrix} a_1,\dots,a_r\\ b_1,\dots,b_s\end{matrix}; 
z\right]=\sum _{k=0} ^{\infty}\frac {\po{a_1}{k}\cdots\po{a_r}{k}}
{k!\,\po{b_1}{k}\cdots\po{b_s}{k}} z^k\ ,
\end{equation*}
we have
\begin{multline*}
\sum_{k} 
{\binom{2n + \ell - k - 4}{\ell -k} } {\binom{n-2}{k}}^2  \\=
\frac {({ \textstyle 2n-3}) _{\ell}} {\ell!}
{{{} _{3} F _{2} \!\left [ \begin{matrix} { 2 - n, 2 - n, -\ell}\\ { 1, 4 - \ell -
      2 n}\end{matrix} ; {\displaystyle 1}\right ] 
      }}={\binom{ \ell + n - 2}{\ell} }^2
\end{multline*}
by means of the Pfaff--Saalsch\"utz summation 
(see \cite[(2.3.1.3),\break Appendix~(III.2)]{Sl})
\begin{equation*} 
{} _{3} F _{2} \!\left [ \begin{matrix} { a, b, -N}\\ { c, 1 + a + b - c -
   N}\end{matrix} ; {\displaystyle 1}\right ]  = 
  {\frac{({ \textstyle c-a}) _{N} \,({ \textstyle c-b}) _{N} } 
    {({ \textstyle c}) _{N} \,({ \textstyle c-a - b}) _{N} }},
\end{equation*}
where $N$ is a nonnegative integer, and similarly
\begin{align*} 
\sum_{k} 
{\binom{2n + \ell - k - 4}{\ell -k} } &
 {\binom{n-3}{k-1}} {\binom{n-1}{k+1}} \\
&=
\frac {({ \textstyle n-2}) _{2} \,({ \textstyle 2n-3}) _{\ell-1} } 
{2\,(\ell-1)!}
{{{} _{3} F _{2} \!\left [ \begin{matrix} { 3 - n, 1 - \ell, 3 - n}\\ { 3, 5 - \ell -
      2\,n}\end{matrix} ; {\displaystyle 1}\right ] 
     }}\\
&={\binom{ \ell + n - 2}{\ell - 1} } {\binom{ \ell + n - 2}{\ell + 1}
},
\end{align*}
again by means of the Pfaff--Saalsch\"utz summation.
%

An alternative, and perhaps more elementary way to prove the equivalence of
\eqref{Eq:Ho1} and \eqref{Eq:Ho2} is to proceed as follows. First, express
both the terms in \eqref{Eq:Ho1} as well as \eqref{Eq:Ho2} in the form
$\sum a_j {\binom{ \ell + 2n - 4 - j}{2n - 4 - j} }$ (where the coefficients
$a_j$ are independent of $\ell$), using for example, Lemmas 3.3 and 3.5 of
\cite{GhorAA}. Then use elementary properties of binomial coefficients such
as those listed in Lemmas 3.1 and 3.2 of \cite{GhorAA} to check that the
corresponding coefficients $a_j$'s are equal.

\medskip

Our next theorem recovers the known formula for the (Krull) dimension 
of a Pfaffian ring (cf\@. \cite{KlLaAA}), and gives a 
closed form expression 
for the multiplicity of a Pfaffian ring. In geometric terms, 
this theorem gives the dimension
(after subtracting $1$) and the degree of
the projective variety in 
$\PP^{{\binom{n}{2}}-1}$ defined by a Pfaffian ideal. 

\begin{theorem}
\label{MainThm2}
The dimension of the Pfaffian ring $R_{r+1}(X)=K[X]/I_{r+1}(X)$  equals 
$r(2n-2r-1)$ and its 
multiplicity $e(R_{r+1}(X))$ equals
\begin{equation} \label{e4.5} 
\prod _{1\le i\le j\le n-2r-1} ^{}\frac {2r+i+j} {i+j}.
\end{equation}
\end{theorem}
\begin{proof} It is well-known that, if the Hilbert series 
of a finitely generated graded 
$K$-algebra $R$ is written
in the form $Q(z)/(1-z)^d$, where $Q(z)$
is a polynomial with rational coefficients
such that $Q(1)\ne 0$, 
then the dimension of $R$ equals $d$ and the multiplicity of $R$ equals $Q(1)$ (see
e.g\@. \cite[Prop.~4.1.9]{BrHeAB}).
(Equivalently, the multiplicity is the sum of the components of the
$h$-vector, the latter being, by definition, the vector of coefficients of
$Q(z)$.) Using Chu--Vandermonde summation 
(see e.g\@. \cite[Sec.~5.1, (5.27)]{GrKPAA}) again, the 
numerator on the right-hand side of \eqref{e4.1} specialized at $z=1$ is
\begin{equation} \label{e4.6} 
{\det\limits_{1\le i,j\le r}\bigg(
\binom {2n-4r}{n-2r-i+j}-
\binom {2n-4r}{n-2r-i-j+1}\bigg)}.
\end{equation}
This determinant can be evaluated by using e.g\@. Theorem~30, (3.18), in
\cite{KratBN}, with $n$ replaced by $r$, 
$q=1$, $A=2n-4r$, and $L_i=-n+2r+i$. The result
is
$$\prod _{i=1} ^{r}\frac {(2n-4r+2i-2)!}
{(n-r-i)!\,(n-r+i-1)!}
\prod _{1\le i<j\le r}
^{}(j-i)\prod _{1\le i\le j\le r}
^{}(i+j-1).
$$
This expression can be transformed into the one given in \eqref{e4.5}.
\quad \quad \qed
\end{proof}

\medskip
\noindent
{\it Example}. 
{}From the Structure Theorem for Gorenstein ideals of codimension $3$ (see, for example
\cite[Sec. 3.4]{BrHeAB}), we know that $I_r(X)$ is {\em the} generic Gorenstein ideal of 
codimension $3$ if $n=2r+1$. In this important special case, the multiplicity formula
of Theorem~\ref{MainThm2} 
simply reduces to
$$
\prod_{1\le i\le j\le 2} \frac{2r - 2 + i+j}{i+j} = \frac{2r (2r+1)(2r+2)}{2\cdot 3 \cdot 4}
=  \frac{r(r+1)(2r+1)}{6}.
$$
This formula is also derived in \cite[p. 29]{HeTrAA} by means of
direct combinatorial considerations. Yet another proof can be found
in \cite{HeTVAA}.

\medskip
\noindent
{\it Remarks}. (1) From the formulae \eqref{e4.1} and \eqref{e4.3}, and the 
observations in the proof of Theorem~\ref{MainThm2}, it is readily seen that the $F_k$'s 
defined in Corollary~\ref{Cor1} give the $h$-vector of the Pfaffian ring $R_{r+1}(X)$.

(2) As we detail in the appendix, a closed form expression for the
multiplicity of a Pfaffian ring had been obtained earlier by Harris
and Tu \cite[Prop.~12]{HT}. Although the form of their expression
(see \eqref{HTProduct}) is different, it is of course completely
equivalent. However, the method with the help of which Harris and Tu
derive their formula is entirely different from ours.

(3) A determinant very similar to the one in \eqref{e4.6}
appears already in \cite[Theorem~5.6]{HeTrAA}
(see \eqref{HerzTrungDet}). However, Herzog and
Trung did not notice that it actually simplifies. (In fact, they
raised the question as to whether or not it simplifies.)

(4) Determinants such as the one in \eqref{e4.6} and the one in
\cite[Theorem~5.6]{HeTrAA} arise quite frequently in the literature, in
particular in connection with the associated counting problem, the
problem of enumerating all families 
$(P_1,P_2,\break\dots,P_r)$ of nonintersecting lattice
paths, where $P_i$ runs from $A_i=(r+i-1,r-i+1)$ to
$E_i=(n-r+i-1,n-r-i+1)$ and does not pass above the line $x=y$, 
$i=1,2,\dots,r$, which we encountered here, or
equivalent problems. For example, if (in contrast to the proof of
\eqref{e4.1a} and
\eqref{e4.2}) we prepend $2(i-1)$ horizontal
steps and append $2(i-1)$ vertical steps to $P_i$, then, out of a
former family of nonintersecting paths, we obtain a family
$(P''_1,P''_2,\dots,P''_r)$ of nonintersecting lattice
paths, where $P''_i$ runs from $A''_i=(r-i+1,r-i+1)$ to
$E''_i=(n-r+i-1,n-r+i-1)$ and does not pass above $x=y$, 
$i=1,2,\dots,r$. See Figure~5 for the
corresponding path family which is obtained out of the one in
Figure~3.b. Again, the number of the
latter families is exactly the same as the number of the former,
because the prepended and appended portions are ``forced," i.e., if 
$(P''_1,P''_2,\dots,P''_r)$ are nonintersecting, they must
contain these prepended and appended portions.

\begin{figure}[h]
$$
\Einheit.4cm
\Koordinatenachsen(12,12)(0,0)
\Kreis(0,0)
\GPunkt(1,1)
\GPunkt(2,1)
\GPunkt(3,1)
\GPunkt(4,1)
\GPunkt(5,1)
\GPunkt(6,1)
\GPunkt(7,1)
\GPunkt(8,1)
\GPunkt(9,1)
\GPunkt(10,1)
\GPunkt(11,1)
\GPunkt(2,2)
\GPunkt(3,2)
\GPunkt(4,2)
\GPunkt(5,2)
\GPunkt(6,2)
\GPunkt(7,2)
\GPunkt(8,2)
\GPunkt(9,2)
\GPunkt(10,2)
\GPunkt(11,2)
\GPunkt(3,3)
\GPunkt(4,3)
\GPunkt(5,3)
\GPunkt(6,3)
\GPunkt(7,3)
\GPunkt(8,3)
\GPunkt(9,3)
\GPunkt(10,3)
\GPunkt(11,3)
\GPunkt(4,4)
\GPunkt(5,4)
\GPunkt(6,4)
\GPunkt(7,4)
\GPunkt(8,4)
\GPunkt(9,4)
\GPunkt(10,4)
\GPunkt(11,4)
\GPunkt(5,5)
\GPunkt(6,5)
\GPunkt(7,5)
\GPunkt(8,5)
\GPunkt(9,5)
\GPunkt(10,5)
\GPunkt(11,5)
\GPunkt(6,6)
\GPunkt(7,6)
\GPunkt(8,6)
\GPunkt(9,6)
\GPunkt(10,6)
\GPunkt(11,6)
\GPunkt(7,7)
\GPunkt(8,7)
\GPunkt(9,7)
\GPunkt(10,7)
\GPunkt(11,7)
\GPunkt(8,8)
\GPunkt(9,8)
\GPunkt(10,8)
\GPunkt(11,8)
\GPunkt(9,9)
\GPunkt(10,9)
\GPunkt(11,9)
\GPunkt(10,10)
\GPunkt(11,10)
\GPunkt(11,11)
\Diagonale(-1,-1){13}
\Pfad(3,3),121121122122\endPfad
\Pfad(2,2),1111211122212222\endPfad
\Pfad(1,1),11111112112122222222\endPfad
\Kreis(3,3)
\Kreis(2,2)
\Kreis(1,1)
\Kreis(9,9)
\Kreis(10,10)
\Kreis(11,11)
\Label\lo{A''_1\kern4pt}(3,3)
\Label\lo{A''_2\kern4pt}(2,2)
\Label\lo{A''_3\kern4pt}(1,1)
\Label\lo{E''_1\kern4pt}(9,9)
\Label\lo{E''_2\kern4pt}(10,10)
\Label\lo{E''_3\kern4pt}(11,11)
\Label\ro{x=y\kern5pt}(12,12)
\hskip4.8cm
$$
\caption{}
\end{figure}

In \cite{DeViAB}, Desainte--Catherine and Viennot showed 
that this counting problem is
equivalent to the problem of counting tableaux with a bounded number
of columns, all rows being of even length. 
(It is from there, that we ``borrowed" the nice
product expression \eqref{e4.5}.)
The counting problem is solved in \cite{DeViAB} (see also
\cite{ChGoAA}) by applying the
main theorem on nonintersecting lattice paths 
\cite[Lemma~1]{LindAA}, \cite[Cor.~2]{GeViAB} and thus obtaining a
determinant, namely the Hankel determinant
$\det_{1\le i,j\le r}(C_{n-2r+i+j-2})$
for the number in question. Here, $C_n$ 
is the $n$-th {\it Catalan number} $\frac {1} {n+1}\binom {2n}n$.
Desainte--Catherine and Viennot evaluate this determinant by means of
the quotient-difference algorithm. However, there are much easier
ways to do it, for example, by using the fact that
$C_n=(-1)^n2^{2n+1}\binom {1/2}{n+1}$ and noticing that therefore
Theorem~26, (3.12), in \cite{KratBN}
is applicable. This latter observation shows that
actually a more general determinant, 
$\det_{1\le i,j\le r}(C_{\lambda_i+j})$, can be evaluated. This
determinant appears also in connection with tableaux counting, see
\cite[second half of Sec.~9]{GeViAB}.
A weighted version of the tableaux counting problem of
Desainte--Catherine and Viennot was solved 
by D\'es\-ar\-m\'enien \cite[Th\'eor\`eme~1.2]{DesaAB}.

\bigskip
\noindent
{\it Acknowledgment}. 
We are grateful to Prof. Piotr Pragacz for 
helpful correspondence.
The first author would like to take this opportunity to gratefully
  acknowledge the support of the Austrian Science Foundation FWF and 
  the warm hospitality of the Institut f\"ur Mathematik der Universit\"at
  Wien for his visit in April--May 1996, which led to this paper.

\appendix

\section*{Appendix: Geometry of Degeneracy Loci and a 
Plethora of Multiplicity Formulae}

In this appendix we briefly review the geometry of degeneracy loci,
and how it leads to multiplicity formulas for the rings that we are 
interested in. For an in-depth treatment of the state-of-the-art in
this subject, we refer the reader 
to \cite{FultPrag}, \cite{Manivel}, and the references therein.

Let us assume, for simplicity, that the ground field is $\C$. 
If $V$ is an $N$-dimensional 
nonsingular projective variety over $\C$, then over the reals, $V$ is a 
compact orientable $(2N)$-dimensional manifold. The top (integral) homology
group $H_{2N}(V)$ is isomorphic to $\Z$, with a generator denoted by
$[V]$ and called the {\it fundamental class} of $V$.
If $W$ is an irreducible closed $d$-dimensional subvariety of $V$, then
the inclusion $i: W \hookrightarrow V$ induces the map 
$$
i_* : \overline{H}_{2d}(W) \to \overline{H}_{2d}(V) = {H}_{2d}(V) 
\simeq {H}^{2e}(V) \ \ {\rm where} \ \ e = \codim W = N-d.
$$
The image of the generator of $\overline{H}_{2d}(W)$ in ${H}^{2e}(V)$ 
is called the {\it fundamental class of $W$ in $V$} and is denoted by $[W]$.
Here, the bar over $H$ indicates that the homology is suitably 
adjusted\footnote{More precisely, instead of singular homology, one uses 
the Borel-Moore homology. 
For details concerning the latter, see \cite[Appendix B]{FultonYT}.}
so that subvarieties $W$ that are singular can also be considered. 
In case $V= \PP^N $, then each $H^{2e}(\PP^N)$ is isomorphic to $\Z$
with a generator given by the fundamental class of a $d$-dimensional
linear subspace or equivalently, by the power\footnote{$H^*(V) = \oplus 
{H}^{i}(V)$ is a graded ring with respect to the cup product, and thus
taking powers makes sense.} $\xi^e$, where $\xi$ denotes the generator
of ${H}^{2}(\PP^N)$ given by the class of a hyperplane. Moreover, in
this case we have 
\begin{equation}
\label{DegClass}
[W] =( \deg W ) \xi^e .
\end{equation}
In this way, the degrees of subvarieties of $\PP^N$ are related to 
their fundamental classes. 

Let $E$ and $F$ be vector spaces over $\C$ of dimensions $n$ and $m$
respectively. An $m\times n$ matrix $A$ over $\C$ defines a map 
$\phi_A : E \to F$ (given by $\phi_A(v) = Av$). For any positive integer
$r \le \min \{m, n\}$, the set 
$$
{\widetilde{M}}_r (m,n) = \{ A \in \C^{mn} : \rank (\phi_A) \le  r \}
$$
is the algebraic variety given by the vanishing of the $(r+1)\times (r+1)$
minors of a generic $m\times n$ matrix. Since the rank condition is
unaltered by changing the matrix $A$ to a nonzero multiple, we may also
consider the corresponding projective variety
\begin{equation}
\label{Mrmn}
{M}_r (m,n) = \{ A \in \PP_C^{mn-1} : \rank (\phi_A) \le  r \}.
\end{equation}
Note that a symmetric $n \times n$ matrix $A$ correspond to maps
$\psi_A : E \to E^*$ where $E^*$ is the dual vector space, satisfying,
$\langle \psi_A x, y\rangle
= \langle \psi_A y, x\rangle$ for all $x,y \in E$, where 
$\langle \ , \ \rangle$ is a dual pairing between $E$ and $E^*$; 
skew-symmetric $n \times n$
matrices have a similar interpretation. The symmetric
or the skew-symmetric matrices will give rise, in the same way, 
to varieties
defined by the minors of a generic symmetric matrix or by the Pfaffians of a 
generic skew-symmetric matrix. 

Geometers like to consider, more generally, the situation when $E$ and $F$ 
are vector bundles, of ranks $n$ and $m$ respectively, over a nonsingular
projective variety (or a complex manifold) $V$, and 
$\phi : E \to F$ is a homomorphism
of vector bundles. Locally, around a point of $V$, the map $\phi$ looks like
$\phi_A$ above. An analogue of (\ref{Mrmn}) is  the set
\begin{equation}
\label{Drphi}
{D}_r (\phi) = \{ x \in V  : \rank (\phi_x) \le  r \}.
\end{equation}
This is a subvariety of $V$, called the {\em degeneracy locus} of rank $r$
associated to $\phi$. There is a natural notion of the dual vector bundle $E^*$
and of symmetric and skew-symmetric bundle homomorphisms $\psi : E \to E^*$;
for such maps, the degeneracy loci are similarly defined and also denoted 
by $D_r(\psi)$. 

It turns out that $D_r(\phi)$ has codimension $\le (m-r)(n-r)$, and when
equality holds (which happens in the generic case), then the
fundamental class $[D_r(\phi)]$ 
is independent of $\phi$ and depends only on
$E$, $F$ and $r$. In fact, it is a polynomial in the so called Chern 
classes\footnote{See \cite[Sec.~3.5.1]{Manivel} for a quick review of Chern
classes and \cite[Sec.~3.2]{FultonIT} for more details.} of $E$ and $F$, 
and is explicitly given by the Giambelli--Thom--Porteous formula:
\begin{equation}
\label{GTP}
[D_r(\phi)] = \det_{1\le i, j\le n-r}  
\left( c_{m-r+j-i}(F-E)  \right).
\end{equation}
Here, $c_{i}(F-E)$ denotes the term of degree
$i$ in the formal expansion of $c(F-E):= c(F)/c(E)$ 
while $c(E)$ denotes the total 
Chern class of $E$, namely, $1+c_1(E)+c_2(E)+ \cdots $.
In particular, if we take $V = \PP_{\C}^{mn-1}$, $E = \O^n_V$, $F = \O_V (1)^m$
and $\phi : E \to F$ the map given by matrix multiplication\footnote{Here, we 
use the standard notation $\O = \O_V$ for the trivial line bundle over $V$, and 
$\O (1)$ for the hyperplane line bundle on $V$; $\O^n$ denotes the direct sum
of $n$ copies of $\O$. Homogeneous coordinates of elements of 
$V = \PP_{\C}^{mn-1}$ can be written as a matrix $A = (a_{ij})$ and then
$\phi (A, v)= Av$.}, then $D_r(\phi)$
is the same as $M_r(m,n)$. In this case $c(E)= 1$ and 
$c_i(F) = {\binom {m} {i}}\xi^i$. Also, $\codim M_r(m,n) = (m-r)(n-r)$. 
Hence from (\ref{DegClass}) and (\ref{GTP}),
we get 
$$
\deg M_r(m,n) = \det_{1\le i, j\le n-r}  
\left(  {\binom{m} {m-r+j-i}} \right) = \prod_{i=0}^{n-r-1}  
\frac{(m+i)! i!}{(r+i)! (r-k+i)!}. 
$$
For a proof of the last equality as well as a proof of (\ref{GTP}), see
\cite[Ch. 2, $\S$ 4]{ACGH}. It may be remarked that equivalent formulae
for the degree of $M_r(m,n)$ also follow from the work of 
Abhyankar \cite{AbhyAB} or from the work of Herzog and Trung 
\cite{HeTrAA}.

The case of varieties defined by the Pfaffians of a skew-symmetric matrix
is similar but with a few twists. 
For a skew-symmetric bundle map $\psi: E \to E^*$, where $E$ 
is as before,
the fundamental class of the degeneracy locus 
$D_{2r}(\psi)$ is given by
the following formula analogous to (\ref{GTP}):
\begin{equation}
\label{HTChern}
[D_{2r}(\psi)] = \det_{1\le i, j\le n-2r-1}  
\left( c_{n-2r+j-2i}(E^*)  \right),
\end{equation}
provided $\psi$ is generic (which basically means that 
$D_{2r}(\psi)$ has
the expected codimension $(n-2r)(n-2r-1)/2$).
An equivalent formula, in terms of the so called Segre 
classes\footnote{\label{f6}%
Formally, the Segre classes $s_i(E)$ of the bundle $E$ can be defined by 
the relation
$\left( \sum_{i=0}^{\infty} s_i(E) t^i \right)
\left( \sum_{i=0}^{\infty} c_i(E) t^i \right) = 1$. Thus $s_0(E)=c_0(E)=1$,
$s_1(E) = -c_1(E)$, $s_2(E) = c_1(E)^2 - c_2(E)$, etc. 
See \cite[Ch.~3]{FultonIT} for more on these.} 
is given by 
\begin{equation}
\label{JLPSegre}
[D_{2r}(\psi)] = \det_{1\le i, j\le n-2r-1}  
\left( s_{2j-i}(E^*)  \right).
\end{equation}
These formulae are due to J\'ozefiak, Lascoux and Pragacz, and
independently, Harris and Tu. To be exact, (\ref{HTChern}) is
given in \cite[Thm.~8]{HT} while  (\ref{JLPSegre}) is a consequence of
\cite[Prop. 5]{JLP}. We can specialize these as before, except that one
has to be a little careful because 
when 
$V= \PP_{\C}^{N-1}$, where $N = n(n-1)/2$, and  $E=\O^n$, the natural
bundle map given by skew-symmetric matrices (whose entries may be viewed
as homogeneous coordinates of elements of $V$) is a skew-symmetric map 
$\psi: \O^n \to (\O^n)^*\otimes \O(1)$. Thus a twist by a line bundle is 
involved\footnote{The appearance of the twist, and the 
corresponding specialization
of Chern classes is explained in the symmetric case in \cite[p. 79]{HT}. 
This depends, in turn, on the `squaring principle' given in 
\cite[pp. 76--78]{HT}. It may be pertinent to remark here that most 
concepts and results discussed so far
extend readily from the complex case to that of an arbitrary ground field 
$K$, at least when $\Char K =0$, if
instead of cohomology rings, one works in the Chow ring of 
algebraic cycles modulo rational equivalence. But it is not clear to us how the
proof of the squaring principle in \cite{HT} would go through in the 
general case. The remarks
in \cite[Sec. 6.4]{FultPrag} may, however, be useful in this context.}.
So in this case one has to 
replace $c_i(E^*)$ by $\frac{1}{2^i} {\binom{n}{i}} \xi^i$ 
and $s_i(E^*)$ by $\frac{1}{2^i} {\binom{n+i-1}{i}} \xi^i$.  
By combining \eqref{DegClass} with \eqref{HTChern} and \eqref{JLPSegre},
this then leads to the following determinantal 
formulae for the 
multiplicity\footnote{It may be noted that the degree of a projective 
variety $V$ is the same as the multiplicity of the cone over $V$ at 
its vertex. Thus the terms {\em degree} and {\it multiplicity} are
sometimes used interchangeably.} of 
the Pfaffian ring $R_{r+1}(X)$.
\begin{eqnarray}
\label{HTDet}
e( R_{r+1}(X)) & =& \frac{1}{2^{\binom{n-2r}{2}}} \det_{1\le i, j\le n-2r-1}  
\left( {\binom{n}{n-2r+j-2i}} \right) \\
& = &\frac{1}{2^{\binom{n-2r}{2}}} \det_{1\le i, j\le n-2r-1}  
\left( {\binom{n+2j-i-1}{2j-i}} \right).
\label{JLPDet}
\end{eqnarray}
The determinant in the first formula has been evaluated by Harris and Tu
\cite[Prop. 12]{HT}, and in this way they get
\begin{equation}
\label{HTProduct}
e( R_{r+1}(X))  = \frac{1}{2^{n-2r-1}} 
\prod_{i = 0}^{n-2r-2}
          \frac{ {\binom{n + i}{2r + 2 i + 1}} }{ {\binom{2i + 1}{i }} }.
\end{equation}
It may be interesting, to compare these with the determinantal formula of
Herzog and Trung \cite{HeTrAA}:
\begin{equation}
\label{HerzTrungDet}
e( R_{r+1}(X))  = 
\det_{1\le i, j\le r}  
\left( {\binom{2n-4r+2}{n-2r-i+j+1}} - {\binom{2n-4r+2}{n-2r-i-j+1}} \right)
\end{equation}
and the product expression given by Theorem \ref{MainThm2}:
\begin{equation}
\label{GKProduct}
e( R_{r+1}(X))  = 
\prod _{1\le i\le j\le n-2r-1} ^{}\frac {2r+i+j} {i+j}.
\end{equation}
Note also that as a consequence of the expressions for the Hilbert 
series in Theorem \ref{MainThm1}, we get three more determinantal 
formulae (e.g., (\ref{e4.6})), which are somewhat similar 
to (\ref{HerzTrungDet}) 
(one of which is given in \eqref{e4.6}). 
Thus one has
indeed a plethora of multiplicity formulae and it may be an interesting
combinatorial exercise to check the resulting identity between any two
by direct methods. To this end, we note that the equivalence of 
(\ref{HTDet}) and (\ref{JLPDet}) 
follows from the well-known Jacobi--Trudi
identities.\footnote{In fact, already the equivalence of \eqref{HTChern} 
and \eqref{JLPSegre} is a consequence of the equivalence 
of the Jacobi--Trudi identities for Schur functions (see e.g\@., 
\cite[Ch.~I, (3.4), (3.5)]{MacdAC} for information on this topic;
see also \cite[Sec.~3.2]{FultPrag} for a geometric point of view). 
The determinants on the right-hand sides of 
(\ref{HTDet}) and (\ref{JLPDet}) can be seen as Schur functions of the
shape $(n-2r-1,n-2r-2,\dots,1)$ because 
the Chern (resp: Segre) classes can be seen as elementary (resp: complete 
homogeneous) symmetric functions in the so called Chern roots; see 
footnote~\ref{f6} and \cite[Ch.~3]{FultonIT} for more details.}
The equivalence of (\ref{HTProduct}) and (\ref{GKProduct}) is
not very difficult to establish directly, 
as also the equivalence
of \eqref{HerzTrungDet} and the similar formulae obtained from
Theorem~\ref{MainThm1}, such as \eqref{e4.6}. But the equivalence of 
(\ref{HTDet}) or (\ref{JLPDet}) with (\ref{HerzTrungDet}) or the similar
formulae obtained from Theorem~\ref{MainThm1} does seem rather intriguing.

\clearpage
\addcontentsline{toc}{section}{Index}
\flushbottom
\printindex

\end{document}